\documentclass[12pt,a4paper]{amsart}

\usepackage{amsfonts,amssymb}
\usepackage{amstext}
\usepackage{amsthm}
\usepackage{amsmath}
\usepackage{amscd}
\newtheorem{thm}{Theorem}[section]
\newtheorem{lem}[thm]{Lemma}
\theoremstyle{definition}

\newtheorem{ex}{exercice}
\newcommand{\be}{\begin{ex} \normal }\newcommand{\ee}{\end{ex}}

\numberwithin{equation}{section}

\newcommand{\N}{\mathbb N}
\newcommand{\R}{\mathbb R}

\newcommand{\ds}{\displaystyle}
\newcommand{\U}{\mathbb U}

\title{Stabilizers of closed sets in the Urysohn space}

\author{Julien Melleray}
\date{}

\begin{document}
 \maketitle

\begin{abstract}
\noindent Building on earlier work of Kat\v{e}tov, Uspenskij
proved in \cite{Uspenskij2} that the group of isometries of
Urysohn's universal metric space $\U $, endowed with the product
topology, is a universal Polish group (i.e it contains an
isomorphic copy of any Polish group). Answering a question of Gao
and Kechris, we prove here the following, more precise result: for
any Polish group $G$, there exists a closed subset $F$ of $\U$
such that $G$ is topologically isomorphic to the group of
isometries of $\U$ which map $F$ onto itself.
\end{abstract}

\begin{section}{Introduction}
\noindent In a posthumously published article (\cite{Urysohn}),
P.S Urysohn constructed a complete separable metric space $\U$
that is \textit{universal} (meaning that it contains an isometric
copy of every complete separable metric space), and
$\omega$-\textit{homogeneous} (i.e such that its isometry group
acts transitively on isometric r-tuples contained in it).\\
In recent years, interest in the properties of $\U$ has greatly
increased, especially since V.V Uspenskij, building on earlier
work of Kat\v{e}tov, proved in \cite{Uspenskij1} that the isometry
group of $\U$ (endowed with the product topology) is a universal
Polish group, that is to say any Polish group is isomorphic to a
(necessarily closed) subgroup of it. \\
In \cite{gaokec}, S. Gao and A.S Kechris used properties of $\U$
to study the complexity of the equivalence relation of isometry
between certain classes of Polish metric spaces; as a side-product
of their construction, they proved the beautiful fact that any
Polish group is (topologically) isomorphic to the isometry group
of some Polish space. A consequence of their construction is that,
for any Polish group $G$, there exists a sequence $(X_n)$ of
closed subsets of $\U$ such that $G$ is isomorphic to $Iso
(\U,(X_n))=\{\varphi \in Iso(\U) \colon \forall n \,
(\varphi(X_n)=X_n)\}$. This led
them to ask the following question (cf \cite{gaokec}): \\
Can every Polish group be represented, up to isomorphism, by a
group of the form $Iso(\U,F)$ for a single subset $F \subseteq \U$
?  \pagebreak

\noindent The purpose of this article is to provide a positive
answer to this question by proving the following theorem:
\begin{thm} \label{answer}
Let $G$ be a Polish group. There exists a closed set $F
\subseteq\U$ such that $G$ is (topologically) isomorphic to
$Iso(F)$, and every isometry of $F$ is the restriction of a unique
isometry of $\U$; in particular, $G$ is isomorphic to $Iso(\U,F)$.
\end{thm}
\noindent This gives a somewhat concrete realization of any Polish
group as a subgroup of $Iso(\U)$.\\
The construction, which will be detailed in section \ref{preuve},
starts with a bounded Polish metric space $X$ such that $G$ is
isomorphic to $Iso(X)$ (the isometry group of $X$, endowed with
the product topology) (Gao and Kechris proved that such an $X$
always exist see \cite{gaokec}). Identifying $G$ with $Iso(X)$, we
construct an embedding of $X$ in $\U$ and a discrete, unbounded
sequence $(x_n) \subseteq \U$ such that $F=X \cup\{x_n\}$ has the
desired properties (here we identify $X$ with its image via the
embedding
provided by our construction).\\

\noindent \textit{Acknowledgements.} Several conversations with
Mathieu Florence while I was working on this paper were very
helpful; for this I am extremely grateful, and owe him many
thanks.

\end{section}

\begin{section}{Notations and definitions} \label{intro}
\noindent If $(X,d)$ is a complete separable metric space, we say
that it is a \textit{Polish metric space}, and often
write it simply $X$. \\
To avoid confusions, we say, if $(X,d)$ and $(X',d')$ are two
metric spaces, that $f$ is an \textit{isometric map} if
$d(x,y)=d'(f(x),f(y))$ for all $x,y \in X$; if $f$ is onto, then
we say that $f$ is an \textit{isometry}.
\\
A \textit{Polish group} is a topological group whose topology is
Polish. If $X$ is a separable metric space, then we denote its
isometry group by $Iso(X)$, and endow it with the product
topology, which turns it into a second countable topological
group, and into a Polish group if $X$ is Polish (see \cite{beckec}
or \cite{Kechris1} for a thorough introduction to the theory of Polish groups). \\
We say that a metric space $X$ is \textit{finitely injective} iff
for any finite subsets $K \subseteq L$ and any isometric map
$\varphi \colon K \to X$ there exists an isometric map $\tilde
\varphi \colon L \to X$ such that $\tilde \varphi_{|_K}=\varphi$.
Up to isometry, $\U$ is the only finitely injective Polish metric
space
(see \cite{Urysohn}).\\
\noindent Let $(X,d)$ be a metric space; we say that $f: X \to \R
$ is a \textit{Kat\v{e}tov map } iff
$$\forall x,y \in X \ |f(x)-f(y)|\leq d(x,y) \leq f(x)+f(y)\ \ .$$
These maps correspond to one-point metric extensions of $X$. We
denote by $E(X)$ the set of all Kat\v{e}tov maps on $X$ and endow
it with the sup-metric, which turns it into a complete
metric space.\\
That definition was introduced by Kat$\check {\mbox{e}}$tov in
\cite{Katetov}, and it turns out to be pertinent to the study of
finitely injective spaces, since one can easily see by induction
that a non-empty metric space $X$ is finitely injective if, and
only if,
$$\forall x_1,\ldots,x_n \in X \ \forall f \in E(\{x_1,\ldots,x_n\})\ \exists z \in X\ \forall i=1\ldots n \ d(z,x_i)=f(x_i)\ .$$
If $Y \subseteq X$ and $f \in E(Y)$, define $\hat{f}\colon X \to
\R$ ( the Kat\v{e}tov extension of $f$) by
$$\hat{f}(x)= \inf\{f(y)+d(x,y) \colon y \in Y\} .$$
Then $\hat{f}$ is the greatest 1-Lipschitz map on $X$ which is
equal to $f$ on $Y$; one checks easily (see for instance
\cite{Katetov}) that $\hat{f} \in E(X) $ and $f\mapsto \hat{f}$ is
an isometric embedding of $E(Y)$ into $E(X)$.\\
To simplify future definitions, we will say that if $f \in E(X)$ and $S \subseteq X$ are such that \\
$\forall x \in X \ f(x)=\inf\{f(s)+d(x,s) \colon s \in S\}$, then
$S$ is a \textit{support of }$f$, or that $S$ \textit{controls} $f$. \\
Notice that if $S$ controls $f\in E(X)$ and $S \subseteq T$, then $T$ controls $f$.\\
Also, $X$ isometrically embeds in $E(X)$ via the Kuratowski map $x
\mapsto f_x$, where \\
$f_x(y)=d(x,y)$. A crucial fact for our purposes is that
$$ \forall f \in E(X)\ \forall x \in X\  d(f,f_x)=f(x).$$
Thus, if one identifies $X$ with its image in $E(X)$ via the
Kuratowski map, $E(X)$ is a metric space containing $X$ and such
that all one-point metric extensions of $X$ embed isometrically in
$E(X)$.\\

\noindent We now go on to sketching Kat\v{e}tov's construction of
$\U$; we refer the reader to \cite{gaokec}, \cite{Gromov},
\cite{Urysohn} and \cite{Uspenskij2} for a more detailed
presentation and proofs of the results we will use below.\\
Most important for the construction is the following
\begin{thm}\label{completion}(Urysohn)
If $X$ is a finitely injective metric space, then the completion
of $X$ is also finitely injective.
\end{thm}
\noindent Since $\U$ is, up to isometry, the unique finitely
injective Polish metric space, this proves that the completion of
any separable finitely injective metric space is isometric to
$\U$.\\
The basic idea of Kat\v{e}tov's construction works like this: if
one lets $X_0=X$, $X_{i+1}=E(X_i)$ then, identifying each $X_i$ to
a subset of $X_{i+1}$ via the Kuratowski map,
let $Y$ be the inductive limit of the sequence $X_i$. \\
The definition of $Y$ makes it clear that $Y$ is finitely
injective, since any $\{x_1,\ldots,x_n\}\subseteq Y$ must be
contained in some $X_m$, so that for any $f\in
E(\{x_1,\ldots,x_n\})$ there exists $z \in X_{m+1}$ such that
$d(z,x_i) =f(x_i)$ for all $i$.\\
Thus, if $Y$ were separable, its completion would be isometric to
$\U$, and
one would have obtained an isometric embedding of $X$ into $\U$.\\
The problem is that $E(X)$ is in general not separable: at each
step, we have added too many functions.\\
Define then $\ds{E(X,\omega)=\{f \in E(X)\colon f \mbox{ is
controlled by some finite } S\subseteq X\}\ .}$ \\
$E(X,\omega)$ is easily seen to be separable if $X$ is, and the
Kuratowski map actually maps $X$
into $E(X,\omega)$, since each $f_x$ is controlled by $\{x\}$.\\
Notice also that, if $\{x_1,\ldots,x_n\} \subseteq X$ and $f \in
E(\{x_1,\ldots,x_n\})$, then its Kat\v{e}tov extension $ \hat{f}$
is in $E(X,\omega)$, and $d(\hat{f},f_{x_i})=f(x_i)$ for all
$i$.\\
Thus, if one defines this time $X_0=X$, $X_{i+1}=E(X_i,\omega)$,
then the inductive limit $Y$ of $\cup X_i$ is separable and
finitely injective, hence its completion $Z$ is
isometric to $\U$, and $X \subseteq Z$. \\
The most interesting property of this construction is that it
enables one to keep track of the isometries of $X$: indeed, any
$\varphi\in Iso(X)$ is the restriction of a unique isometry
$\tilde \varphi$ of $E(X,\omega)$, and the mapping $\varphi
\mapsto \tilde \varphi$ from $Iso(X)$ into $Iso(E(X,\omega))$ is a
continuous group
embedding.\\
That way, we obtain for all $i$ continuous embeddings
$\Psi^i\colon Iso(X) \to Iso(X_i)$, such that
$\Psi^{i+1}(\varphi)_{|_{X_i}}=\Psi^i(\varphi)$ for all $i$
and all $\varphi \in Iso(X)$.\\
This in turns defines a continuous embedding from $Iso(X)$ into
$Iso(Y)$, and since extension of isometries defines a continuous
embedding from the group of isometry of any metric space into that
of its completion (see \cite{Uspenskij1}), we actually have a
continous embedding of $Iso(X)$ into the isometry group of $Z$,
that is to say $Iso(\U)$
(and the image of any $\varphi \in Iso(X)$ is actually an extension of $\varphi$ to $Z$ ).\\

\end{section}
\begin{section}{Proof of the main theorem} \label{preuve}
\noindent To prove theorem \ref{answer} ,we will use ideas very
similar to those used in \cite{gaokec} ;
all the notations are the same as in section \ref{intro}.\\
We will need an additional definition, which was introduced in \cite{gaokec}:\\
If $X$ is a metric space and $i \geq 1$, let
$$E(X,i)=\{f \in E(X): f \mbox{ has a support of cardinality } \leq i\}$$
We endow $E(X,i)$ with the sup-metric. \\
Gao and Kechris proved the following result, of which we will give
a new, slightly simpler proof:
\begin{thm}\label{gaokec} (Gao-Kechris) \\
If $X$ is a Polish metric space and $i \geq 1$ then $E(X,i)$ is a
Polish metric space.
\end{thm}
$ $ \\
\noindent {\bf Proof:} \\
Notice first that the separability of $E(X,i)$ is easy to prove;
we will prove  its completeness by induction on $i$.\\
The proof for $i=1$ is the same as in \cite{gaokec}; we include it for completeness.\\
First, let $(f_n)$ be a Cauchy sequence in $E(X,1)$. \\
It has to converge uniformly to some Kat\v{e}tov map $f$, and it
is enough to prove that $f \in
E(X,1)$.\\
By definition of $E(X,1)$, there exists a sequence $(y_n)$ such
that $$ \forall x \in X \ f_n(x)= f_n(y_n)+d(y_n,x) \qquad \qquad
(*)$$ But then let $\varepsilon >0$, and let $M$ be big
enough that $m,n \geq M \Rightarrow d(f_n,f_m) \leq \varepsilon$. \\
Then, for $m,n \geq M$, one has
$$2d(y_n,y_m)=(f_n(y_m)-f_m(y_m))+(f_m(y_n)-f_n(y_n)) \leq 2 \varepsilon\ .$$
This proves that $(y_n)$ is Cauchy, hence has a limit $y$.\\
One easily checks that $f(y)= \lim f_n(y_n)$, so that $(*)$ gives
us, letting $n \to \infty$
$$ \forall x \in X \ f(x)=f(y)+d(y,x)\qquad .$$

$ $\\
\noindent That does the trick for $i=1$; suppose now we have
proved the result for $1 \ldots i-1$, and let $(f_n)$ be a Cauchy
sequence in $E(X,i)$.\\
By definition, there are $y^n_1, \ldots y^n_i$ such that:\\
 $\forall x \in X\  f_n(x)=\ds{\min_{1\leq j \leq
i}\{f_n(y^n_j)+d(y^n_j,x) \}}\qquad (**)$ . \\

\noindent Once again, $(f_n)$ converges uniformly to some
Kat\v{e}tov map $f$,
and we want to prove that $f \in E(X,i)$. \\
Thanks to the induction hypothesis, we can assume that there is
$\delta >0$ such that for all $n$ and all $k \neq j\leq i$ one has
$d(y^n_j,y^n_k) \geq 2 \delta$ (if not, a subsequence of
$(f_n)$ can be approximated by  a Cauchy sequence in $E(X,i-1)$, and the induction hypothesis applies).\\
Let $d_n=\min\{f_n(x) \colon x \in X\}$.\\
Then $(d_n)$ is Cauchy, so it has a limit $d \geq 0$; up to some
extraction, and if necessary changing the enumeration of the
sequence,  we can assume that there is $p \geq 1$ and $\delta' >0$ such that:\\
- $\forall j \leq p \ f_n(y^n_j) \to d$\\
- $\forall j >p\ \forall n \ f_n(y^n_j)> d+ \delta'$. \\
Let $\varepsilon >0$, $\alpha=\min(\delta,\delta',\varepsilon) $
and choose $M$ big enough
that $n,m \geq M \Rightarrow d(f_n,f_m) < \frac{\alpha}{4}$ and $|f_n(y^n_j)-d| < \frac{\alpha}{4}$ for all $j \leq p$.\\
Then, for $n,m \geq M$ and $j \leq p$ one has:\\
$f_n(y^m_j)<d+\frac{\alpha}{2}$, so there exists $k \leq p$
such that $f_n(y^m_j)=f_n(y^n_k)+d(y^m_j,y^n_k)$. \\
Such a $y^n_k$ has to be at a distance strictly smaller than
$\delta$ from $y^m_j$: there is at most one $y^n_k$ that can work,
and there is necessarily one. Thus, one obtains, as in the case
$i=1$, that $d(y^n_k,y^m_j) \leq  \varepsilon $ .\\
This means that one can assume, choosing an appropriate
enumeration, that for $k\leq p$ each sequence $(y^n_i)_n$ is
Cauchy, hence has a limit $y_k$ .\\
Define then $\ds{\tilde{f}_n \colon x \mapsto \min_{1\leq k\leq p}
\{f_n(y^n_i)+d(x,y^n_k)\}}$ .
\\
$\tilde{f}_n \in E(X,p)$, and one checks easily, since $y^n_k \to
y_k$ for all $k \leq p$, that $(\tilde{f}_n)$ converges uniformly
to $\tilde{f}$, where $\ds{\tilde{f} (x)= \min_{1\leq k\leq p}
\{f(y_k)+d(x,y_k)\} \ .}$ \\

\noindent If $p=i$ then we are finished; otherwise, notice that,
using again the induction hypothesis, we may assume that there is
$\eta
>0$ such that
$$\forall n \, \forall j>p \qquad \ds{f_n(y^n_j) <\tilde{f}_n(y^n_j)-\eta \qquad (***)}.$$
Now define $\tilde{g}_n$ by $\ds{\tilde{g}_n(x)= \min_{j>p}
\{f_n(y^n_j)+d(x,y^n_j)\}}$. \\
Choose $M$ such that $\ds{n,m \geq M \Rightarrow d(f_n,f_m)<
\frac{\eta}{4}}$ and $\ds{d(\tilde{f}_n,\tilde{f}_m) <\frac{\eta}{4}}$. \\

\noindent Then (***) shows that for, all $n,m \geq M$ and all
$j>p$,
$$f_m(y^n_j) \leq f_n(y^n_j)+\frac{\eta}{4} \leq \tilde{f}_n(y^n_j)-\frac{3\eta}{4}
\leq \tilde{f}_m(y^n_j)-\frac{\eta}{2}\ ,$$ so that
$f_m(y^n_j)=f_m(y^m_k) + d(y^n_j,y^m_k) $ for some $k >p$.
\\
Consequently, for $m,n \geq M$ and $j >p$,
$f_m(y^n_j)=\tilde{g}_m(y^n_j) $; by definition, $f_m(y^m_j)=\tilde{g}_m(y^m_j) $.\\
This proves that for all $n,m \geq M$ one has
$d(\tilde{g}_n,\tilde{g}_m) \leq d(f_n,f_m)$, so that
$(\tilde{g}_n)$ is Cauchy in $E(X,i-p)$, hence has a limit
$\tilde{g} \in E(X,i-p)$ by the induction hypothesis.\\
But then, (**) shows that, for all $x \in X$, $\ds{f(x)= \min
(\tilde{f}(x),\tilde{g}(x))}$, and this concludes the proof.
$\hfill
\lozenge$\\

\noindent If $Y $ is a nonempty, closed and bounded subset of a
metric space $X$, define
$$E(X,Y)=\{f \in E(X) \colon \exists d \in \R^{+}\ \forall x \in X
\ f(x)=d+d(x,Y)\}$$ \noindent $E(X,Y)$ is closed in $E(X)$, and is
isometric to a closed, unbounded interval of $\R^+$.\\

\noindent Now we can go on to the\\
\noindent {\bf Proof of theorem \ref{answer}.} \\
Essential to our proof is the fact that for every Polish group $G$
there exists a Polish space $(X,d)$ such that $G$ is isomorphic to
the group of isometries of $X$ (This result was proved by Gao and Kechris, see \cite{gaokec}). \\
So, let $G$ be a Polish group, and $X$ be a metric space
such that $G$ is isomorphic to $Iso(X)$.\\
One can assume that $X$ contains more than two points, and $(X,d)$
is bounded, of diameter $d_0 \leq 1$.(If not, define
$\ds{d'(x,y)=\frac{d(x,y)}{1+d(x,y)}}$. Then $(X,d')$ is now a
bounded Polish metric space with the same topology as $X$, and the
isometries of $(X,d')$ are exactly the isometries of
$(X,d)$ ).\\

\noindent Let $X_0=X$, and define inductively bounded Polish
metric spaces $X_i$, of diameter $d_i$, by:\\
$$X_{i+1}=\big\{f\in E(X_i,i)\cup \bigcup_{j<i} \ E(X_i,X_j)\colon \forall x \in X_i \ f(x)\leq 2 d_i\big\}$$ (We endow
$X_{i+1}$ with the sup-metric; since $X_i$ canonically embeds
isometrically in $X_{i+1}$ via the
Kuratowski map, we assume that $X_i \subseteq X_{i+1}$).\\
Remark that $d_i \to + \infty$ with $i$, and that each $X_i$ is a
Polish metric space. \\
Let then $Y$ be the completion of $\ds{\bigcup_{i \geq 0} X_i}$. \\
The definition of  $\cup X_i$ makes it easy to see that it is
finitely injective, so that $Y$ is
isometric to $\U$.\\
Also, any isometry $g \in G$ extends to an isometry of $X_i$, and
for any $i$ and $g \in G$ there is a unique isometry $g^i$ of
$X_i$ such that $g^i(X_j)=X_j $ for all $j \leq i$ and
${g^{i}}_{|_{X_0}}=g$ (same proof as in  \cite{Katetov}). \\
Remark also that the mappings $g \mapsto
g^i$, from $G$ to $Iso(X_i)$, are continuous (see \cite{Uspenskij1}).\\
All this enables us to assign to each $g$ an isometry $g^{*}$ of
$Y$, given by ${g^*}_{|_{X_i}}=g^i$, and this defines a continuous
embedding of $G$ into $Iso(Y)$ (see again \cite{Uspenskij1} for details). \\
It is important to remark here that, if $f \in X_{i+1}$ is defined
by \\
$f(x)=d+d(x,X_j)$ for some $d \geq 0$ and some $j < i$, then
$g^*(f)=f$ for all $g \in G$ (This was the aim of the definition
of $X_i$: adding "many" points that are fixed by the action of $G$).\\

\noindent Notice that an isometry $\varphi$ of $Y$ is equal to
$g^{*}$ for some $g\in G$ if, and only if, $
\varphi(X_n)=X_n$ for all $n$.\\
The idea of the construction is then simply to construct a closed
set $F$ such that $\varphi(F)=F$ if, and only if,
$\varphi(X_n)=X_n$ for all $n$. To achieve this, we will build $F$
as a set of carefully chosen "witnesses".\\

\noindent The construction  proceeds as follows:\\
First, let $(k_i)_{i \geq 1}$ be an enumeration of the
non-negative integers where every number appears infinitely many times.\\
Using the definition of the sets $X_i$, we choose recursively for
all $i \geq 1$ points $a_i \in \cup_{n \geq 1}X_n$(the
witnesses), non-negative reals $e_i$, and a nondecreasing sequence of integers $(j_i)$ such that :\\

\noindent- $e_1 \geq 4; \  \forall i \geq 1\ e_{i+1}> 4 e_i $ . \\
- $\forall i \geq 1 \ j_i \geq k_i \, , \  \ a_i \in X_{j_i+1}
\mbox{ and } \forall x \in X_{j_i} \ d(a_i,x)=e_i+d(x,X_{k_i-1})$ \\
- $\forall i \geq 1 \ \forall g \in G \  g^{*}(a_i)=a_i$.\\
(This is possible, since at step $i$ it is enough to fix $e_i
>\max(4e_{i-1},\mbox{diam}(X_{k_i}))$, then find $j_i\geq
\max(1+j_{i-1},k_i)$ such that $\mbox{diam}(X_{j_i}) \geq  e_i$,
and define $a_i \in X_{j_i+1}$ by the equation above; then, by
definition of $g^{*}$ and of $a_i$, one has
$g^{*}(a_i)=a_i$ for all $g \in G$)\\

\noindent Let now $\ds{F= X_0
\cup \{a_i\}_{i\geq 1}}$; since $X_0$ is complete, and $d(a_i, X_0)=e_i \to + \infty$, $F$ is closed. \\
We claim that for all $\varphi \in Iso(Y)$, one has
$$(\varphi(F)=F) \iff (\varphi \in G^{*}).$$
The definition of $F$ makes one implication obvious.\\
To prove the converse, we need a lemma:

\begin{lem} \label{sympa}
If $\varphi \in Iso(F)$, then $\varphi(X_0)=X_0$, so that
$\varphi(a_i)=a_i$ for all $i$. Moreover, there exists $g \in G$
such that $\varphi={g^{*}}_{|_F}$.
\end{lem}

\noindent Admitting this lemma for a moment, it is now easy to
conclude: \\
Notice that lemma \ref{sympa} implies that $G$ is isomorphic to
the isometry group of $F$, and that any isometry of
$F$ extends to $Y$.\\
Thus, to conclude the proof of theorem \ref{answer}, we only need
to show that the extension of a given isometry of $F$ to $Y$ is
unique. As explained before, it is enough to show that, if
$\varphi \in Iso(Y)$ is such that $\varphi(F)=F$, then
$\varphi(X_n)=X_n$ for all $n \geq 0$. \\
So, let $\varphi \in Iso(Y)$ be such that $\varphi(F)=F$. \\
It is enough to prove that $\varphi(X_n) \supseteq X_n$ for all $n
\in \N$ (since this will also be true for $\varphi^{-1}$), so
assume that this is not true, i.e there is some $n \in \N$ and $x
\not \in X_n$ such that $\varphi(x)\in X_n$. \\
Let $\delta=d(x,X_n)>0$ (since $X_n$ is complete), and pick $y \in
\cup X_m$ such that $d(x,y)\leq \frac{\delta}{4}$.\\
Then $y \in X_m \setminus X_n$ for some $m > n$; now choose $i$
such that $k_i=n+1$ and $j_i \geq m$ .\\
Then we know that\\
$\ds{d(\varphi(y),\varphi(a_i))=d(y,a_i)=e_i+d(y,X_n) \geq
e_i+\frac{3 \delta}{4}}$, and \\
$\ds{d(a_i,\varphi(y))\leq d(a_i,\varphi(x))+ d(x,y) \leq
e_i+\frac{\delta}{4}}$, so that
$d(\varphi(a_i),a_i) \geq \frac{\delta}{2}$, and this contradicts lemma \ref{sympa}. $ \hfill \lozenge$ \\

\noindent It only remains to give the \\
{\bf Proof of lemma \ref{sympa}: }\\
Since we assumed that $X_0$ has more than two points and
$\mbox{diam}(X_0) \leq 1$, the definition of $F$ makes it clear
that
$$\forall x \in F \ (x \in X_0) \Leftrightarrow (\exists y \in F \colon \ 0 < d(x,y) \leq 1 )$$
The right part of the equivalence is invariant by isometries of
$F$, so this proves that $\varphi(X_0)=X_0$ for any $\varphi \in
Iso(F)$. In turn, this easily implies that
$\varphi(a_i)=a_i$ for all $i \geq 1$. \\
Thus, if one lets $g \in G$ be such that
$g_{|_{X_0}}=\varphi_{|_{X_0}}$, we have shown that
$\varphi={g^*}_{|_F}$. $\hfill \lozenge$

\end{section}

\noindent  Equipe d'Analyse Fonctionelle, Universit\'e Paris 6 \\
Boîte 186, 4 Place Jussieu, Paris Cedex 05, France.\\
e-mail: melleray@math.jussieu.fr

\begin{thebibliography}{7}
\bibitem{beckec} H.Becker and A.S Kechris, \textit{The Descriptive Set Theory of Polish Group
Actions}, London Math. Soc. Lecture Notes Series, 232, Cambridge
University Press (1996).
\bibitem{gaokec} S. Gao and A.S Kechris, \textit{ On the
classification of Polish metric spaces up to isometry}, Memoirs of
Amer. Math. Soc., 766, Amer. Math. Soc. (2003).
\bibitem{Gromov} M. Gromov, \textit{Metric Structures for Riemannian and Non-Riemannian
spaces.} Birkauser (1998)
\bibitem{Katetov} M. Kat\v{e}tov, \textit{ On universal metric spaces}, Proc. of the 6th Prague Topological Symposium (1986), Frolik
(ed). Helderman Verlag Berlin, pp 323-330 (1988).
\bibitem{Kechris1} A.S Kechris, \textit{Classical descriptive set
theory}, Springer-Verlag (1995).
\bibitem{Kechris3} A.S Kechris, \textit{Actions of Polish groups and Classification
Problems}, Analysis and Logic, London Math. Soc. Lecture Notes
Series,Cambridge University Press (2000).
\bibitem{Urysohn} P.S Urysohn, \textit{ Sur un espace m\'etrique
universel}, Bull. Sci. Math 51 (1927), pp 43-64 and 74-96.
\bibitem{Uspenskij1} V.V Uspenskij, \textit{Compactifications of topological
groups}, Proc. of the 9th Prague Topological Symposium (2001) pp
331-346, Topology Atlas, Toronto (2002).
\bibitem{Uspenskij2} V.V Uspenskij, \textit{On the group of isometries of the Urysohn universal metric
space}, Comment. Math. Univ. Carolinae, 31(1) (1990).
\end{thebibliography}
\end{document}